\documentclass{amsart}
\usepackage{amsmath}

\theoremstyle{plain}
\newtheorem{thm}{Theorem}[section]
\newtheorem{prop}[thm]{Proposition}
\newtheorem{lemma}[thm]{Lemma}
\newtheorem{cor}[thm]{Corollary}

\theoremstyle{definition}

\newcommand{\PP}{{\mathbb P}}
\newcommand{\F}{{\mathbb F}}
\newcommand{\PPP}{{\mathbb P}^3}
\newcommand{\Ff}{{\mathbb F}_8}
\newcommand{\pr}{\PPP_{\Ff}}
\newcommand{\pgl}{{\rm PGL}_4(\F_8)}
\newcommand{\Aut}{{\rm Fix}}
\newcommand{\GL}{{\rm GL}}

\newcommand{\Z}{{\mathbb Z}}

\begin{document}

\begin{center}

\textbf{THE MAXIMUM NUMBER OF POINTS ON A CURVE OF GENUS $4$ OVER $\F_8$ 
IS $25$}

\vskip 0.2cm

\textsc{David Savitt}\footnote{Partially supported by an NSERC postdoctoral
fellowship.}

\vskip 0.2cm

with an Appendix by \textsc{Kristin Lauter}
\end{center}

\vskip 1.0 cm

\markboth{\textsc{David Savitt}}{\textsc{The maximum number of points on a 
curve of genus $4$ over $\F_8$ is $25$}}
\pagestyle{myheadings}

\section{Introduction}

Our aim in this paper is to prove that a smooth geometrically 
irreducible curve $C$ of genus $4$
over the finite field $\F_8$ may have at most $25$ $\F_8$-points.  Our
strategy is as follows: if $C$ has more than $18$ $\F_8$-points, then 
$C$ may not be hyperelliptic, and so the canonical divisor of $C$ yields 
an embedding
of $C$ into $\PPP_{\F_8}$.  The image of $C$ under this embedding is
a degree $6$ curve which is precisely the intersection of an irreducible 
cubic hypersurface with an 
irreducible quadric
hypersurface, both defined over $\F_8$.  (This is Example IV.5.2.2 in 
\cite{hartshorne}.  Hartshorne works over an algebraically closed field, 
but his argument is equally valid over the smaller field.  See, for 
example, Theorem III.5.1 in \cite{hartshorne} and Theorem A.4.2.1 in 
\cite{sh} for the necessary tools.)

Consequently, finding the maximum possible number of points on a curve of 
genus $4$ over $\F_8$ is reduced to a finite task: one can write 
down all cubic hypersurfaces and all quadric hypersurfaces in $\PPP$, and 
count the 
number of points on their intersection.  As a practical matter, however, 
one must make significant reductions before this program becomes 
computationally feasible.  For example, the space of homogeneous cubics 
in four variables is already $\binom{6}{3} - 1 = 19$-dimensional.  

We begin in section \ref{quadrics} by noting that up to isomorphism
there are only three irreducible quadric surfaces in $\PPP_{\F_8}$ which
contain many $\F_8$-points.  Therefore we may select representatives of 
the isomorphism classes and assume that our curve $C$
lies on one of these three specific quadrics.  Next, we recall (see 
\cite{lauter} and
\cite{gv}) that it is known that any curve of genus $4$ over $\F_8$ has no
more than $27$ points, and that such curves with $25$ points exist.  
Moreover, using the techniques of \cite{lauter2}, K. Lauter demonstrates
in an appendix to this paper that such curves with $26$ points do not 
exist. We may therefore 
suppose that the curve $C$ for which we are searching has exactly 
$27$ points.  In section \ref{reductions}, we employ the following 
strategy to reduce the problem further.  If $Q$ is one of our three 
quadrics, then the subgroup $\Aut(Q)$ of $\pgl$ preserving $Q$ is large.  
If $P$ is a cubic surface and if $\sigma \in \Aut(Q)$, then $ P \cap Q$ 
and $\sigma(P) \cap Q = \sigma(P \cap Q) $ have the 
same number of points.  If the intersection $P \cap Q$ is a geometrically 
irreducible curve of degree $6$, then by B\'ezout's theorem the 
intersection may contain at most three points of any line. We study the 
action of $\Aut(Q)$ on the points of $Q$ to show that if $S \subset Q$ is 
a subset with $27$ points, no four of which are colinear,  
then we may find $\sigma \in \Aut(Q)$ such that $\sigma(S)$ 
contains a particular list of points of $Q$ (or one of several lists of 
points of $Q$).

The problem is therefore reduced to studying cubics $P$ 
which contain particular points of $Q$, cutting down significantly on the 
dimension of the space of cubics under consideration.  Depending on the  
cubic, we are able to eliminate between $5$ and $7$ dimensions in this
fashion.  The space is cut down further by $4$ dimensions by noting that 
we may subtract appropriate multiples of our quadric $Q$.  Thus we have 
reduced a $19$-dimensional search space over $\F_8$ to a 
search space over $\F_8$ of no greater than $10$ dimensions, which is 
easily tractable for a computer.

Finally, we note that this search will a priori turn up many cubics and
quadrics whose intersection contains $27$ points.  This is
because we will find many reducible (or at least geometrically reducible)  
intersections.  These ``bad'' curves are relatively straightforward to 
identify and discard. In section \ref{bad}, we give a precise list of the 
ways in which bad curves with $27$ points can occur.

\section*{Acknowledgements}

The author is grateful to J-P. Serre for his comments and corrections, and
in particular for the suggestion that section \ref{bad} be included.  We also
thank Jason Starr for several helpful conversations,  William Stein for 
the use
of his computer, and the anonymous referee for
his or her comments.  Computations were performed partly by C programs, and
partly using the MAGMA package.  This problem came to the author's
attention at the 2000 Arizona Winter School on Arithmetic Algebraic
Geometry, and the author thanks the organizers of this conference for
their hard work and hospitality.

\section{Quadric surfaces in $\pr$} \label{quadrics}

Let $C$ be a non-hyperelliptic curve of genus $4$ over $\F_8$.  As
we have noted, we may suppose that $C$ is canonically embedded into
$\pr$ as the intersection of an irreducible quadric hypersurface $Q$
with an irreducible cubic hypersurface $P$.  It is a classical result 
that over a finite field $\F$, there are exactly three 
reduced and geometrically 
irreducible quadric surfaces in $\PPP_{\F}$ up to $\F$-isomorphism: the 
split nonsingular quadric (isomorphic to $\PP^1_{\F} \times \PP^1_{\F}$),
the nonsplit nonsingular quadric (the quadratic twist of $\PP^1_{\F}
\times \PP^1_{\F}$), and the singular quadric.

We give an argument, essentially found on p. 206 of 
\cite{ACGH}, explaining for each $C$ into which of the above categories 
the quadric $Q$ falls.  Note that any linear system of degree $3$ and 
dimension at least $1$ on
$C$ defines a ruling of $Q$.  Indeed, if $D$ is a divisor in such a linear
system, then by the geometric version of the Riemann-Roch theorem, the
linear span in $\PPP$ of the support of $D$ is a line.  By B\'ezout's theorem,
this line is contained in $Q$.

The $\F_8$-scheme $W^{1}_{3}(C)$ defined in \cite{ACGH}, 
whose geometric points correspond
to the complete linear series of degree $3$ and dimension at least $1$
on $C$, is a zero-dimensional affine scheme, and by the Thom-Porteous 
formula this scheme has degree $2$.  Hence
there are exactly three possibilities for $W^{1}_{3}(C)$: two reduced
$\F_8$-points (so $Q$ is the split nonsingular quadric), two conjugate 
$\F_{64}$-points (nonsplit nonsingular), and one nonreduced
$\F_8$-point (singular).

To make our classification of quadrics concrete, we first recall the 
following result from \cite{arf}:

\begin{prop}\label{woof}  Let $\F$ be a field of characteristic $2$.  Then any 
quadratic form in $n$ variables over $\F$ is equivalent to one of the form
$$ \sum_{i=1}^{\mu} x_i y_i + \sum_{j=\mu+1}^{\mu+\nu} (a_j x_j^2 + x_j 
y_j + b_j y_j^2) + \sum_{k=1}^{d} c_k z_k^2 $$ with $2\mu + 2\nu + d \le 
n$.
\end{prop}

This is by no means a 
classification: two distinct quadratic forms written as above may 
still be isomorphic.  For example, when the field $\F$ is perfect
evidently we 
may take $d=0$ or $1$ and $c_1=1$.  Similarly we may 
suppose each $a_j=1$.

When the field $\F=\F_{2^n}$ with $n$ odd, one can check 
with little difficulty that the form $x^2 + xy + by^2$ is equivalent
either to the form $xy$ or to $x^2 + xy + y^2$, depending on whether
or not the form nontrivially represents $0$ over $\F$.  Combining
this with the identity
$$X^2 + XY + Y^2 + Z^2 = XY + (X+Y+Z)^2$$
and the fact that
$$ (X^2 + XY + Y^2) + (Z^2 + ZW + W^2)$$
is identically equal to
$$(X+Z+W)(Y+Z+W) + (X+Y+Z)(X+Y+W) \,, $$ we obtain the following
version of Proposition \ref{woof}.

\begin{prop} \label{eight} Let $\F_{2^n}$ be the finite field with $2^n$ elements with 
$n$ an odd integer.  Then any quadratic form over $\F_{2^n}$ 
is equivalent over $\F_{2^n}$ to a form with one of the following shapes:
\begin{itemize}
\item $\sum_{i=1}^{\mu} x_i y_i$
\item $\sum_{i=1}^{\mu} x_i y_i + (X^2 + XY + Y^2)$
\item $\sum_{i=1}^{\mu} x_i y_i + (Z^2)$.
\end{itemize}
\end{prop}

We are interested in particular in the geometrically integral quadric 
surfaces in 
$\pr$, which correspond to geometrically irreducible quadratic forms in at 
most four variables over $\F_8$.  Their classification is as follows.

\begin{prop} \label{three} Up to $\F_8$-isomorphism, there are exactly 
three 
geometrically irreducible 
quadratic forms 
in four variables $X,Y,Z,W$ over $\F_8$.  They are: $XY+ZW$ (the
split non-degenerate form), $X^2 + XY + Y^2 + ZW$ (the non-split
non-degenerate form), and $XY+Z^2$ (the degenerate form).
\end{prop}

\begin{proof} According to Proposition \ref{eight}, 
up to isomorphism there are at 
most six quadratic forms in four variables over any finite field 
$\F_{2^n}$ with $n$ odd, namely: $XY$, $XY + ZW$, $X^2 + XY + Y^2$, $X^2 
+ XY + Y^2 + ZW$, $Z^2$, and $Z^2 + XY$.  The forms $XY$ and $Z^2$ are 
reducible, and $X^2 + XY + Y^2$ is irreducible but geometrically 
reducible, and so we eliminate them.

The hypersurface defined by $XY + ZW$ is $\F_8$-isomorphic
to $\PP^1 \times \PP^1$, and possesses two $\F_8$-rulings.  The 
hypersurface defined by $X^2 + XY + Y^2 + ZW$ is $\F_{64}$-isomorphic
to $\PP^1 \times \PP^1$, and so one sees that it has two 
Galois-conjugate rulings over $\F_{64}$ but contains no lines over $\F_8$.
Finally, the hypersurfaces defined by $XY + ZW$ and $X^2 + XY + Y^2 + ZW$
are non-singular, whereas the hypersurface defined by $XY + Z^2$ is 
singular at $[0:0:0:1]$.  
These facts together show that these three forms cannot be $\F_8$-isomorphic.
\end{proof}

\noindent \textbf{Remark.}  We can also see that these forms are not
$\F_8$-isomorphic by
verifying that a different number of points of $\pr$ lie on each of the
resulting quadric surfaces.  In fact there are $81$ points on the surface
$XY + ZW = 0$, there are $73$ points on the surface $X^2 = YZ$, and there are
$65$ points on the surface $X^2 + XY + Y^2 = ZW$.

\section{Reductions} \label{reductions}

\subsection{Action of $\pgl$ on quadrics}

In this subsection, we describe the subgroups of $\pgl$ preserving each 
of our quadrics.  If we can correctly list these subgroups in their 
entirety, we will automatically be able to obtain a proof that the 
description is correct, by counting the size of the orbits of our quadrics 
under $\pgl$.

We begin with the quadric $XY + ZW = 0$.  This quadric is 
isomorphic to 
$\PP^1 \times \PP^1$, as can be seen via the map $\PP^1 \times \PP^1 
\rightarrow \{XY=ZW\}$ sending $([x:y],[z,w]) \mapsto [xz:yw:xw:yz]$.  The 
inverse map is defined on coordinate patches, for example sending 
$[X:Y:Z:W] \mapsto ([X:W],[W:Y])$ on the affine $\{W \neq 0\}$.  The group 
${\rm PGL}_2(\F_8) \times {\rm PGL}_2(\F_8) \times C_2$ acts on $\PP^1 
\times \PP^1$, where the cyclic factor $C_2$ is generated by 
an automorphism interchanging the two copies of $\PP^1$.  Evidently 
each nontrivial one of these automorphisms yields a nontrivial element of 
$\pgl$ preserving $XY+ZW=0$.

We turn next to the quadric $XY = Z^2$.  One may easily check that the map
$$
\begin{pmatrix}
X \\ 
Y \\ 
Z \\ 
W
\end{pmatrix}
\mapsto
\begin{pmatrix}
a & b & 0 & 0 \\
c & d & 0 & 0 \\
\sqrt{ac} & \sqrt{bd} & \sqrt{ad+bc} & 0 \\
* & * & * & e
\end{pmatrix}
\begin{pmatrix}
X \\ 
Y \\ 
Z \\ 
W
\end{pmatrix}
$$
preserves $XY = Z^2$, where 
$\begin{pmatrix}
a & b \\
c & d
\end{pmatrix}$ 
is an element of ${\GL}_2(\F_8)$, $e \in \F_8^{\times}$, 
and each $* \in \F_8$.  These will be all the elements of $\Aut(XY=Z^2)$.

Next, we verify that $\Aut(X^2 + XY + Y^2 = ZW)$ acts doubly-transitively
on $\F_8$-points of the quadric.  Indeed, we claim that for any point $p$
on $X^2 + XY + Y^2 =ZW$ other than $[0:0:0:1]$, there is an element of
$\Aut(X^2 + XY + Y^2 = ZW)$ sending $[0:0:1:0]$ to $p$ while fixing 
$[0:0:0:1]$.  Then
for any pair of points $p_1,p_2$ we may send $p_1$ to $[0:0:1:0]$, then
use the automorphism interchanging $W$ and $Z$ to map $p_1$ to
$[0:0:0:1]$.  If $p_2$ has now been moved to $p_3$, we finish via a map
preserving $[0:0:0:1]$ and sending $p_3$ to $[0:0:1:0]$, so the pair
$(p_1,p_2)$ has been moved to $([0:0:0:1],[0:0:1:0])$, and the group is
doubly-transitive.

To see the claim, notice that for an element $x \in \F_8$, the map sending 
$X \mapsto X + xZ, Y \mapsto Y, Z \mapsto Z, W \mapsto W + xY + x^2 Z$ 
preserves $X^2 + XY + Y^2 = ZW$, sends $[0:0:1:0]$ to $[x:0:1:x^2]$, and 
fixes $[0:0:0:1]$.  Now the map 
sending $X \mapsto X, Y \mapsto Y + yZ, Z \mapsto Z, W \mapsto W + yX + 
y^2 Z$ preserves $X^2 + XY + Y^2 = ZW$, sends $[x:0:1:x^2]$ to 
$[x:y:1:x^2+xy+y^2]$, and fixes $[0:0:0:1]$.  Since 
$[x:y:1:x^2+xy+y^2]$ is a general point on the curve besides $[0:0:0:1]$, 
this proves the claim.  

Now an element of $\GL_4(\F_8)$ preserving $X^2 + XY + Y^2 = ZW$ and 
fixing $[0:0:0:1]$ and $[0:0:1:0]$ will be of the form
$$\begin{pmatrix}
a & b & 0 & 0 \\
c & d & 0 & 0 \\
0 & 0 & z & 0 \\
0 & 0 & 0 & w
\end{pmatrix}$$
where 
$z,w \in \F_8^{\times}$,
$\begin{pmatrix}
a & b \\
c & d
\end{pmatrix}$ is an element of $\GL_2(\F_8)$ preserving the form $X^2 + 
XY + Y^2 = 0$, and $a,b,c,d,z$ determine $w$.  One checks that there are 
exactly $126$ such elements of $\GL_2(\F_8)$.
They are the scalar multiples of 
the 
following $18$ matrices: the identity matrix, the 
matrix 
$\begin{pmatrix}
0 & 1 \\
1 & 0
\end{pmatrix}$, 
the four matrices with three entries equal to $1$ and the 
other equal to $0$, and, for each of the three roots $\eta$ of $\eta^3 + 
\eta + 1 = 
0$, the four $90$-degree rotations of the matrix
$\begin{pmatrix}
\eta & \eta^2 \\
\eta^{-3} & \eta
\end{pmatrix}$.  

Furthermore, $X^2 + XY + Y^2 = 0$ in 
$\PPP$ has automorphisms given by completing those $126$ matrices 
$\begin{pmatrix}
a & b \\
c & d
\end{pmatrix}$ to matrices
$$\begin{pmatrix}
a & b & 0 & 0 \\
c & d & 0 & 0 \\
* & * & * & * \\
* & * & * & *
\end{pmatrix}$$
where the last two rows are independent of the first two.  

We now verify that we have indeed found all of the automorphisms of these 
quadrics.  

\begin{itemize}

\item For $XY+ZW$, we have found $2 \cdot ((8^2 - 1)(8^2 - 8)/7)^2 = 
508032$ automorphisms.  This has index $68024320$ in $\pgl$, which is 
therefore an upper bound on the size of the orbit of $XY+ZW$ in the space 
of quadric surfaces.  

\item For $XY=Z^2$, we have found 
$(8^2-1)(8^2-8)\cdot 8^3 \cdot 7 / 7= 1806336$ automorphisms, giving an 
upper bound of $19131840$ on the orbit.  

\item For $X^2 + XY + Y^2 = ZW$, we 
have found $65 \cdot 64 \cdot 126 \cdot 7 / 7 = 524160$ automorphisms, 
giving an upper bound of $65931264$ on the orbit.  

\item For $X^2 + XY + Y^2$, 
we have found $126 * (8^4-8^2)*(8^4-8^3) / 7 = 260112384$ automorphisms, 
giving an upper bound of $132860$ on the orbit.  

\item It is easy to 
see that the form $X^2$ has orbit of size $(8^4-1)/7 = 585$ and $XY$ has 
orbit of size $(8^4-1)(8^4-8)/(2\cdot7^2) = 170820$.  
	
\end{itemize}

Finally, we note 
that 
$68024320+19131840+65931264+132860+585+170820=153391689=(8^{10}-1)/7$, 
precisely the number of quadric surfaces, and so we confirm that we have 
indeed found all the automorphisms of these quadrics.

\subsection{Reductions for $XY+ZW$}

Observe that $\{XY+ZW=0\} \cong \PP^1 \times \PP^1$ is a ruled surface, 
and in particular that the set of $\F_8$-points of $\PP^1 \times \PP^1$ 
may be written as the union of the nine lines $\{l\} \times \PP^1_{\F_8}$, 
for $l \in \PP^1_{\F_8}$, and as the union of the nine lines $\PP^1_{\F_8} 
\times \{r\}$ for $r \in \PP^1_{\F_8}$.  Each of these lines on $\PP^1 
\times \PP^1$ maps to a line on $\{XY+ZW=0\}$.  

In the remainder of this subsection, we suppose that a cubic hypersurface
$P \subset \pr$ intersects the quadric $\{XY+ZW=0\}$ in a smooth
geometrically irreducible curve $C$ with $27$ $\F_8$-points.

If $P$ intersected any of these lines on $\{XY+ZW=0\}$ in at least $4$
points, then by B\'ezout's theorem the line would be contained in $P$, and
consequently the line would be contained in the intersection $P \cap
\{XY+ZW=0\}$.  Therefore the curve $C$ would be reducible, which we have 
assumed is not the case.  We may therefore conclude that $P$ intersects 
each of these lines in at most $3$ points.  However, since there are nine 
lines in each ruling, $P$ must intersect each 
of these lines in \textit{exactly} $3$ points.  Note that this 
argument yields 
a 
combinatorial proof that if the canonical embedding of a smooth curve of 
genus $4$ over $\F_8$ lies on $\{XY+ZW=0\}$, then it cannot contain $28$ 
points.

Write the $\F_8$-points of $\PP^1 \times \PP^1$ as $(l_i,r_j)$ with $0 \le
i,j \le 8$.  We have seen that for each $i$ there are exactly three $j$
such that $(l_i,r_j)$ lies on $P$, and similarly for each $j$ there are
exactly three $i$.  Suppose, after renumbering, that $(l_0,r_0)$,
$(l_0,r_1)$, and $(l_0,r_2)$ all lie on $P$.  We divide into two cases.  
First, suppose there exists $i>0$ such that two of
$(l_i,r_0),(l_i,r_1),(l_i,r_2)$ lie on $P$. After renumbering, we may
assume may assume that $(l_i,r_0),(l_i,r_1)$ lie on $P$, and we may select
$i' \neq 0,i$ so that $(l_{i'},r_2)$ lies on $P$.  Since ${\rm
PGL}_2(\F_8)$ acts $3$-transitively on $\PP^1_{\F_8}$, we may select an
automorphism $\sigma$ of $\PP^1 \times \PP^1$ such that
$([0:1]:[0:1]),([0:1]:[1:0]),([0:1]:[1,1]),([1:0],[0:1]),([1:0],[1:0]),
([1:1],[1:1])$ all lie on $\sigma(P)$.  Therefore, without loss of
generality, in this case we may assume that these six points lie on $P$.
We refer to this as the $3,2,1$-case.

Second, suppose that no such $i$ exists.  Without loss of generality, 
after renumbering we may assume that 
$(l_1,r_0),(l_2,r_0),(l_3,r_1),(l_4,r_1),(l_5,r_2),(l_6,r_2)$ all lie on 
$P$.  Then, by the pigeonhole principle, for some $j > 2$ there are $1\le 
i,i' \le 6$ so that $(l_i,r_j)$ and $(l_{i'},r_j)$ lie on $P$.   If 
$\{i,i'\} = \{1,2\},\{3,4\}$, or $\{5,6\}$, we may suppose after 
renumbering that $\{i,i'\}=\{1,2\}$, and we are reduced to the case of the 
previous paragraph: namely 
$(l_0,r_0),(l_1,r_0),(l_2,r_0),(l_1,r_j),(l_2,r_j),$ and some 
$(l_0,r_{j'})$ lie on $P$, so after interchanging the two copies of 
$\PP^1$ and applying an element of ${\rm PGL}_2(\F_8) \times {\rm 
PGL}_2(\F_8)$, we may again assume that 
$([0:1]:[0:1]),([0:1]:[1:0]),([0:1]:[1,1]),([1:0],[0:1]),([1:0],[1:0]),
([1:1],[1:1])$ all lie on $P$.

On the other hand, if $\{i,i'\} \neq \{1,2\},\{3,4\}$, or $\{5,6\}$, we 
may assume (after renumbering) that $\{i,i'\}=\{1,3\}$.  In this case we 
have $(l_0,r_0)$, $(l_0,r_1)$, $(l_1,r_0),$ $(l_1,r_j)$, $(l_3,r_1)$, 
and $(l_3,r_j)$ 
all lying 
on 
$P$.  Applying an element of ${\rm PGL}_2(\F_8) \times {\rm 
PGL}_2(\F_8)$ we may assume that $ ([0:1],[0:1]), ([0:1],[1:0]), 
([1:0],[0:1]), ([1:0],[1:1]), ([1:1],[1:0]), ([1:1],[1:1])$ all lie on 
$P$.  We refer to this as the $2,2,2$-case.
Moreover, we may suppose that $([1:0],[1:0])$ is not on $P$, or else 
we would be able to reduce to the $3,2,1$-case.

Suppose that homogeneous cubic polynomial defining $P$ is written $c_{X^3} 
X^3 + c_{X^2 Y} X^2 Y + \cdots + c_{W^3} W^3$.  We can now verify the 
following proposition.

\begin{prop} \label{red1} If there exists a cubic hypersurface $P \subset 
\pr$ whose 
intersection with $\{XY=ZW\}$ is a smooth geometrically irreducible curve 
of genus $4$ with $27$ $\F_8$-points, then there exists such a 
hypersurface whose coefficients satisfy one or the other of 
the two sets of conditions below:
\begin{enumerate}
\item
\begin{itemize}
\item $c_{X^3}=c_{Y^3}=c_{Z^3}=c_{W^3}=c_{X^2 
Y}=c_{XY^2}=c_{Z^2W}=c_{ZW^2}=0$,
\item $c_{Y^2 W}=c_{YW^2}=1$, and
\item $c_{X^2Z}+c_{X^2W}+c_{XYZ}+c_{XYW}+c_{XZ^2}+c_{XZW}
+c_{XW^2}+c_{Y^2Z}+c_{YZ^2}+ c_{YZW} = 0$, or
\end{itemize}
\item
\begin{itemize}
\item $c_{X^3}=1$,
\item $c_{Y^3}=c_{Z^3}=c_{W^3}=c_{X^2 Y}=c_{XY^2}=c_{Z^2W}=c_{ZW^2}=0$,
\item $c_{XZ^2} = c_{X^2Z}+1$, $c_{XW^2} = c_{X^2W} + 1$, and
\item $c_{XYZ} + c_{XYW} + c_{XZW}+ 
c_{Y^2Z}+c_{Y^2W}+c_{YZ^2}+c_{YZW}+c_{YW^2} = 1$.
\end{itemize}
\end{enumerate}
\end{prop}

\begin{proof}  Recall that we map $\PP^1 \times \PP^1 \rightarrow 
\{XY=ZW\}$ via $([x:y],[z,w]) \mapsto [xz:yw:xw:yz]$. In the 
$3,2,1$-case, we have shown that we may assume  $([0:1]:[0:1]),
([0:1]:[1:0]),([0:1]:[1,1]),([1:0],[0:1]),([1:0],[1:0]),
([1:1],[1:1])$ all lie on $P$.  In $\PPP$-coordinates, these six points 
are, respectively,
$[0:1:0:0]$,$[0:0:0:1]$,$[0:1:0:1]$,$[0:0:1:0]$,$[1:0:0:0]$, 
and $[1:1:1:1]$.
For these points to lie on $P$, it follows that 
$c_{X^3}=c_{Y^3}=c_{Z^3}=c_{W^3}=0$, that $c_{Y^2W}=c_{YW^2}$, and that 
all $20$ coefficients sum to zero.  If $c_{Y^2 W} = c_{YW^2} = 0$, one
easily verifies that the line $[0:Y:0:W]$ is contained in the 
curve, 
and so we may suppose without loss of generality that $c_{Y^2 
W}=c_{YW^2}=1$.  Further, 
by subtracting appropriate 
multiples of the quadric $XY=ZW$, we may suppose that $c_{X^2 Y} = 
c_{XY^2}=c_{Z^2W}=c_{ZW^2}=0$.  

In the $2,2,2$-case we may assume that $ ([0:1],[0:1]), ([0:1],[1:0]), 
([1:0],[0:1]), ([1:0],[1:1]), ([1:1],[1:0]), ([1:1],[1:1])$ all lie on 
$P$.  In $\PPP$-coordinates, these six points are, respectively,
$[0:1:0:0]$,$[0:0:0:1]$,$[0:0:1:0]$,$[1:0:1:0]$,$[1:0:0:1]$, and 
$[1:1:1:1]$.  
For these points to lie on $P$, it follows that 
$c_{Y^3}=c_{W^3}=c_{Z^3}=0$, that $c_{X^3} + c_{X^2Z} + c_{XZ^2} = 0$, 
that 
$c_{X^3}+c_{X^2W}=c_{XW^2}=0$, and that all the coefficients sum to $0$.
Moreover, we may assume that $([1:0],[1,0])$, which in 
$\PPP$-coordinates is $[1:0:0:0]$, does not lie on $P$.  This implies that 
$c_{X^3}\neq 0$, so we may suppose without loss of generality that 
$c_{X^3}=1$.  Once again,
by subtracting appropriate multiples of the quadric $XY=ZW$, we may 
suppose that $c_{X^2 Y} = c_{XY^2}=c_{Z^2W}=c_{ZW^2}=0$.  
\end{proof}

\subsection{Reductions for $XY=Z^2$}

Suppose that a cubic hypersurface $P \subset \pr$ intersects the quadric 
$\{XY=Z^2\}$ in a smooth geometrically irreducible curve $C$ with $27$ 
$\F_8$-points.

The $\F_8$-points of the surface $\{XY=Z^2\}$ are ruled by the pencil of 
nine lines $[1:z^2:z:W]$,$[0:1:0:W]$ parametrized by the variable $W$, all 
passing through the point 
$[0:0:0:1]$.  By an argument essentially the same as the pigeonhole 
argument in the 
previous subsection, we see that $[0:0:0:1]$ cannot lie on $P$, while 
each of the nine lines intersects $C$ in exactly $3$ other $\F_8$-points.
Note that once again we obtain an elementary proof that there cannot be 
$28$ points on such a curve $C$ lying on this quadric.

We remark that the collection of affine transformations of $\F_8$, i.e.,
the set of maps $x \mapsto e x + f$ with $f \in \F_8$,
$e \in \F_8^{\times}$, acts transitively on the set of $3$-element
subsets of $\F_8$.  Notice that there are $56$ affine transformations of 
$\F_8$ and $56$ $3$-element subsets of $\F_8$, so it suffices to prove 
that the stabilizer of the $3$-element subset $\{0,1,\eta\}$ is trivial.  
(Recall that $\eta$ is a chosen root of $\eta^3 + \eta + 1=0$.)  This is 
easy to check.  For example, the affine transformation swapping $0$ and 
$1$ is $x \mapsto 1-x$, which does not fix $\eta$; and the affine 
transformation sending $0$ to $1$ and $1$ to $\eta$ is $x \mapsto 
(\eta-1) x + 1$, which does not send $\eta$ to $0$.  

Now, for each line $l_z = \{[1:z^2:z:W]\}$, let $S_z = \{ W \ | \ [ 
1:z^2:z:W] \in C \}$.  Observe that each $S_{z}$ has size $3$, and so 
there is a unique transformation $x \mapsto e_z x + f_z$ mapping $S_z$ to
$\{0,1,\eta\}$.  Since there are eight $S_{z}$'s, by the pigeonhole 
principle some element $e \in \F_8^{\times}$ occurs twice in the list of 
$e_z$'s.  
Suppose $e = e_{z_1} = e_{z_2}$.  Choose any element 
$\begin{pmatrix}
a & b \\
c & d
\end{pmatrix}$ 
of ${\rm 
GL}_2(\F_8)$ 
sending the vectors $(1,z_1^2),(1,z_2^2)$ to $(0,1),(1,0)$ respectively.
Suppose that 
this matrix maps the 
line $[1:z_3^2]$ to the line $[1:1]$.  (What we say below will 
work equally well in the case that the transformation maps the line 
$[0:1]$ to the line $[1:1]$, which we omit for ease of notation.)
Select any point of the form $[1:z_3^2:z_3:w_3]$ on $C$.
Then we can solve the system of equations
\begin{eqnarray*}
g_X + g_Y z_1^2 + g_Z z_1 & = & f_{z_1} \\
g_X + g_Y z_2^2 + g_Z z_2 & = & f_{z_2} \\
g_X + g_Y z_3^2 + g_Z z_3 & = & e w_3
\end{eqnarray*}
for the variables $g_X,g_Y,g_Z$.  
Let $\sigma$ be the transformation $$
\begin{pmatrix}
X \\ 
Y \\ 
Z \\ 
W
\end{pmatrix}
\mapsto
\begin{pmatrix}
a & b & 0 & 0 \\
c & d & 0 & 0 \\
\sqrt{ac} & \sqrt{bd} & \sqrt{ad+bc} & 0 \\
g_X & g_Y & g_Z & e
\end{pmatrix}
\begin{pmatrix}
X \\ 
Y \\ 
Z \\ 
W
\end{pmatrix} \,.
$$
Then $\sigma$ preserves $XY=Z^2$, and we have constructed $\sigma$ so that 
$\sigma(P)$ contains the seven points 
$[0:1:0:0]$, $[0:1:0:1]$, $[0:1:0:\eta]$, $[1:0:0:0]$, $[1:0:0:1]$,
$[1:0:0:\eta]$, $[1:1:1:0]$.  Then the following proposition holds.

\begin{prop} \label{red2} If there exists a cubic hypersurface $P \subset 
\pr$ whose 
intersection with $\{XY=Z^2\}$ is a smooth geometrically irreducible curve 
of genus $4$ with $27$ $\F_8$-points, then there exists such a 
hypersurface whose coefficients satisfy the 
following conditions:
\begin{itemize}
\item $c_{X^3} = c_{X^2Y} = c_{XY^2} = c_{Y^3} = c_{Z^3} = c_{Z^2W} = 0$
\item $c_{W^3}=1$, $c_{X^2W}=c_{Y^2W}=\eta$, $c_{XW^2}=c_{YW^2}=\eta^3$.
\item $c_{X^2Z}+c_{XYZ}+c_{XZ^2}+c_{Y^2Z}+c_{YZ^2}=0$.
\end{itemize}
\end{prop}

\begin{proof}  We have seen that the under the hypothesis of the 
proposition, there exists such a 
hypersurface $P$
containing the above seven points and not containing the point 
$[0:0:0:1]$.  From the latter, we may assume without loss of generality 
that $c_{W^3}=1$.  Subtracting appropriate multiples of the quadric 
$XY=Z^2$, we may assume $c_{X^2Y} = c_{XY^2} = c_{Z^3} = c_{Z^2W} = 0$.  
Since $[0:1:0:0]$ and $[1:0:0:0]$ are on the cubic $P$, we get 
$c_{X^3}=c_{Y^3}=0$.  From the presence of $[0:1:0:1]$ on the cubic $P$, 
we get $c_{Y^2W} + c_{YW^2} + 1 = 0$.  From the presence of $[0:1:0:\eta]$ 
on the cubic $P$, we get $c_{Y^2W} \eta + c_{YW^2} \eta^2 + \eta^3 = 0$.  
It follows that $c_{Y^2W}=\eta$ and $c_{YW^2}=\eta^3$.  Similarly 
$c_{X^2W}=\eta$, $c_{XW^2}=\eta^2$.  The last condition follows from 
previous deductions and the presence of $[1:1:1:0]$ on the cubic $P$.
\end{proof}

\subsection{Reductions for $X^2 + XY + Y^2 = ZW$}

Suppose that a cubic hypersurface $P \subset \pr$ intersects the quadric 
$\{X^2+XY+Y^2=ZW\}$ in a smooth geometrically irreducible curve $C$ with 
$27$ $\F_8$-points.

Since $\Aut(X^2 + XY + Y^2 = ZW)$ acts $2$-transitively on the points of
$X^2 + XY + Y^2 = ZW$, we may assume without loss of generality that
$[0:0:1:0]$ and $[0:0:0:1]$ lie on $P$.  Recall that the elements of 
$\pgl$
preserving $X^2 + XY + Y^2 = ZW$ and fixing those two points are of the
form 
$$\begin{pmatrix} 
a & b & 0 & 0 \\ 
c & d & 0 & 0 \\ 
0 & 0 & z & 0 \\
0 & 0 & 0 & w \end{pmatrix}$$ 
where 
$\begin{pmatrix} 
a & b \\ 
c & d
\end{pmatrix}$ 
preserves the form $X^2 + XY + Y^2$, and so such elements
of $\pgl$ permute the nine conics $C_y = [1:y:Z:(1+y+y^2)Z^{-1}]$, $y \in 
\F_8$, and
$C_{\infty} = [0:1:Z:Z^{-1}]$, each conic parametrized by the variable
$Z$, and each conic passing through the two points $[0:0:1:0]$ and
$[0:0:0:1]$.  Unfortunately our previous pigeonhole arguments do not seem
to be of value here, because we now would need seven points of one of
these conics to lie on the curve $C$ to induce a contradiction.

One checks, using our explicit list of the $18$ elements of ${\rm 
PGL}_2(\F_8)$ preserving $X^2 + XY + Y^2=0$, that the action of $\Aut(X^2 
+ XY + Y^2=ZW)$ on the set of nine conics is as follows: the subsets 
$\{C_0,C_1,C_{\infty}\}$, $\{C_{\eta},C_{\eta^2},C_{\eta^{-3}}\}$, and 
$\{C_{\eta^{-1}},C_{\eta^{-2}},C_{\eta^{3}}\}$ are always permuted as 
blocks, and the action on the set of three blocks is the cyclic group of 
order 
$3$.  The stabilizer of each of each block induces the full symmetric 
group of order $6$ on the three elements of the block.

By the pigeonhole principle, since there are $25$ points of $C$ (besides 
$[0:0:0:1]$ and $[0:0:1:0]$) on the nine curves, it follows that the 
conics in at least one of the blocks contain a total of at least $9$ 
points of $C$.  Permuting the blocks, we may assume that this block is 
$\{C_0,C_1,C_{\infty}\}$.  Permuting the conics within the block, we may 
also 
assume that $$\#(C_{\infty} \cap P) \ge \#(C_0 \cap P) \ge \#(C_1 \cap 
P) \,. $$  

Certainly we now have $\#(C_{\infty} \cap P) \ge 3$.  Applying 
transformations of the form $X \mapsto X, Y \mapsto Y, Z \mapsto \alpha Z, 
W \mapsto \alpha^{-1} W$ and transformations of the form $X \mapsto X, Y 
\mapsto Y, Z \mapsto \alpha W, W \mapsto \alpha^{-1} Z$, as well as by 
applying the Frobenius element of $\F_8$ to the coefficients of $P$, we 
may suppose that $\#(C_{\infty} \cap P)$ contains the two points 
$[0:1:1:1]$ and $[0:1:\eta:\eta^{-1}]$, and at least one of the two points 
$[0:1:\eta^2:\eta^{-2}]$ and $[0:1:\eta^3:\eta^{-3}]$.  We obtain the 
following proposition.

\begin{prop} \label{red3}  If there exists a cubic hypersurface $P \subset 
\pr$ whose 
intersection with $\{X^2 + XY + Y^2 =ZW\}$ is a smooth geometrically 
irreducible curve 
of genus $4$ with $27$ $\F_8$-points, then there exists such a 
hypersurface satisfying the following conditions:
\begin{itemize}
\item $c_{X^3} = c_{X^2 Y} = c_{X^2 Z} = c_{X^2 W} = c_{Z^3} = c_{W^3} = 
0$,
\item $c_{Y^2Z}= \eta^{-1} c_{Y^2W} + \eta^3 c_{YZ^2} + \eta c_{YW^2} 
 + c_{Z^2W} + \eta^{-1} c_{ZW^2}$,
\item $c_{Y^3} = 
c_{Y^2Z}+c_{Y^2W}+c_{YZ^2}+c_{YZW}+c_{YW^2}+c_{Z^2W}+c_{ZW^2}$,
\item at least one of $[0:1:\eta^2:\eta^{-2}$ and $[0:1:\eta^3:\eta^{-3}]$ 
lies on $P$,
\item $\#(C_{\infty} \cap P) \ge \#(C_0 \cap P) \ge \#(C_1 \cap
P)$ and $\#(C_{\infty} \cap P) + \#(C_0 \cap P) + \#(C_1 \cap
P) \ge 9$.
\end{itemize}
\end{prop}

\begin{proof}
Subtracting appropriate multiples of the quadric, we may assume that 
$c_{X^3}=c_{X^2Y}=c_{X^2Z}=c_{X^2W}=0$.  Since we may assume that 
$[0:0:0:1]$ and $[0:0:1:0]$ lie on the cubic $P$, it follows that we may 
suppose $c_{Z^3}=c_{W^3}=0$.  The two long sums ensure that $[0:1:1:1]$ 
and $[0:1:\eta:\eta^{-1}]$ lie on $P$.  That we may suppose the remainder 
of the conditions follows from our reductions preceeding the proposition.
\end{proof}

\section{Computations}

\subsection{Publicly available data}

The programs we use, the data they produce, and documentation, are 
available on the web at
\begin{center}
{\sf http://www.math.mcgill.ca/$\sim$dsavitt/curves/}
\end{center}
and the longest of our computations took under two days to run.

\subsection{Listing cubics}

The computations we perform are straightforward.  We write a C program 
to 
perform arithmetic in $\F_8$, and then for each of our three quadrics, we 
simply cycle through all possibilities for the coefficients of homogeneous 
cubics in four variables subject to the conditions we are able to impose 
from propositions \ref{red1}, \ref{red2}, and \ref{red3}.  For each 
possible vector
of 
coefficients, we count how many points of the quadric under consideration 
lie on the cubic.  
Each time the intersection 
contains exactly $27$ points, the program 
prints the cubic polynomial in a format which is readable by the MAGMA 
computation package.  In order to speed this up significantly, we store 
in advance the value of each cubic monomial evaluated at each $\F_8$-point 
of the quadric, so that to determine whether a point of the quadric lies 
on the cubic is simply a matter of evaluating a predetermined linear form 
in the coefficients.  For the quadric $X^2 + XY + Y^2 = ZW$, we also add 
routines to check the final two conditions of proposition \ref{red3} and 
discard those cubics in violation of them. 

In order to build redundancy into our computations, we write MAGMA
routines which given a cubic will count the number of points of our
quadric which lie on that cubic.  Using these routines, we can confirm
that our C programs are correctly counting the points on our cubics;  
indeed we can list the points on the cubic and check that the points we
wish to force to lie on the cubic are really there.  However, the 
streamlined C programs will be faster than the MAGMA routines, which is 
why we use the C program and not MAGMA for the computations.

\subsection{Discarding cubics}

From the above computations, we obtain a long list of cubics whose 
$\F_8$-intersection with a particular quadric has size $27$.  If it is 
true that there are no smooth geometrically irreducible curves of 
genus $4$ over $\F_8$ with exactly $27$ points, we expect that each of 
these intersections will be (geometrically) reducible.  In order to test 
this, for each of these cubic-quadric pairs we use MAGMA to count the 
number of $\F_{64}$-points on their intersection.  If the original curve 
were actually smooth and geometrically irreducible, then the number of 
$\F_{64}$-points will be one of the possibilities admitted by the Weil 
conjectures.  If the original curve is reducible, then we expect the 
number of $\F_{64}$-points will be too large.

Explicitly, the methods of section $2$ of \cite{lauter2} leave only two 
possibilities for the list of eigenvalues of Frobenius on a smooth 
geometrically irreducible curve of genus $4$ over $\F_8$ with $27$ 
$\F_8$-points.  If the eigenvalues are $\alpha_i,\overline{\alpha}_i$, 
$i=1,2,3,4$, the possibilities are: $(- \alpha_i - \overline{\alpha}_i)_i 
= (5,5,5,3)$ and $(- \alpha_i - \overline{\alpha}_i)_i = (\frac{9 \pm 
\sqrt{5}}{2} , \frac{9 \pm \sqrt{5}}{2})$.  Using that $\alpha_i 
\overline{\alpha_i} = 8$, we compute that $\sum_i (\alpha_i^2 + 
\overline{\alpha}_i^2) = 20$ or $22$, and so the total number of 
$\F_{64}$-points must be either $1 + 64 - 20 = 45$ or $1 + 64 - 22 = 43$.

In fact, our computations in MAGMA show that every one of the 
cubics we have listed intersects the associated quadric in at least $119$ 
points.  This establishes:

\begin{thm}  There is no smooth, geometrically irreducible curve of genus 
$4$ over $\F_8$ with $27$ points.
\end{thm}

Combined with what was already known, we obtain:

\begin{cor} The maximal number of points on a curve of genus $4$ over 
$\F_8$ is $25$.
\end{cor}

\noindent \textbf{Remark.}  It would be of interest to know whether the 
combinatorial arguments we have given which eliminate the possibility 
of $28$ points on an irreducible curve of genus $4$ over $\F_8$ lying 
on $XY=Z^2$ or $XY=ZW$ can be improved to eliminate the possibility of 
$27$ points, or can be extended to curves lying on $X^2 + XY + Y^2=ZW$.

\section{Bad curves with $27$ points} \label{bad}

As explained above, in our computer search we find numerous examples where
our cubic and our quadric intersect in exactly $27$ $\F_8$-points.  
However, when we count the number of $\F_{64}$-points on the intersection, 
we find that the answer is always in the following list: $119$, $181$, 
$189$, $191$, $195$, $197$, $199$, or $205$.  Moreover, on the 
degenerate and the non-split non-degenerate quadrics, we only find 
examples with $189$ and $191$ $\F_{64}$-points.  In this section, we 
explain 
why these are the only possibilities, and we list (along with 
examples) precisely the ways in which they can occur.  This provides 
significant reassurance that our computer calculations are correct.

\subsection{Preliminary lemmas}

For ease of reference, we note the following facts:

\begin{lemma} \label{field} If $K/k$ is any nontrivial field extension, 
then a curve of 
degree $d$ over $K$ which is not definable over $k$ may have at most $d^2$ 
$k$-points under any embedding into $\PPP_{K}$.
\end{lemma}

\begin{proof} 
By B\'ezout's theorem, two plane curves of degree $d$
intersect in $d^2$ points.  As a consequence, two different curves of
degree $d$ in projective space may intersect in at most $d^2$ points:
otherwise, they coincide under every projection to the plane, and so they
must coincide. As a result, there is at most one curve of degree $d$
through any $d^2+1$ points in projective space.  However, if there is only
one curve of degree $d$ through a set of $k$-points, then by linear
algebra that curve is defined over $k$.  The lemma follows.
\end{proof}

\begin{lemma} \label{max}
If the intersection of a cubic and a quadric in $\pr$ has an 
component defined over $\F_8$ and of degree $3$, $4$, or $5$, then that 
component has at 
most $9$, $14$, or $18$ $\F_8$-points respectively.
\end{lemma}

\begin{proof}
Any component of our intersection which is a plane curve lies on a 
quadric, and so has degree at most $2$.  Therefore any cubic component has 
genus $0$, any quartic component has genus at most $1$, and any quintic 
component has genus at most $2$.  (See Figure 18 on page 354 of 
\cite{hartshorne}.)  The Serre-Weil bounds on the number of points on
curves of genus $0$, $1$, and $2$ over $\F_8$ are $9$, $14$, and $19$
respectively.  The first two of these bounds are met.  The maximum
number of points on a curve of genus $2$ over ${\mathbb F}_q$ was
determined for all $q$ by Serre (this is Th\'eor\`eme 4 in \cite{serre2},
and may also be found as Proposition 1 in \cite{gv}).  When $q=8$, this 
bound is $18$.
\end{proof}

\begin{lemma} \label{red} If the intersection $C$ of a cubic and a 
quadric in $\pr$ 
has $27$ $\F_8$-points but is not a smooth, geometrically irreducible 
curve of 
genus $4$, then the intersection is geometrically reducible.
\end{lemma}

\begin{proof} Assume that $C$ is geometrically irreducible but singular.
We will 
show that it cannot have $27$ points. Since the intersection is not 
planar, the arithmetic genus is at most $4$.  (Again, see Figure 18 in 
\cite{hartshorne}.)  Let $C'$ be the normalization of $C$.  Then by the 
discussion in section IV.7 of \cite{serre}, the 
arithmetic genus of $C'$ is $4-a$ where $a$ is an integer between $1$ and 
$4$, and moreover the number of $\F_8$-points of $C'$ differs from the 
number of $\F_8$ points by at most $a$.  By the Weil conjectures, $C$ may 
have at most $9 + 5\cdot(4-a) + a = 29 - 4a \le 25$ points.
\end{proof}

Similarly, suppose $C$ is a singular curve over $\F_8$ of arithmetic genus 
$1$.  Then the normalization $C'$ has arithmetic genus $0$, so has exactly 
$9$ $\F_8$-points.  The singularity of $C$ must be an ordinary 
double-point, and the number of $\F_8$-points of $C$ must be either $8$ or 
$10$, depending on whether the points of $C'$ lying over the singularity 
are defined over $\F_8$ or $\F_{64}$ respectively.  In either case, the 
number of points of $C$ over $\F_{64}$ will be $64$.

Finally, we note that the components of a geometrically reducible curve 
are permuted by Galois.  In particular, if there is only one component 
of a curve over $\F_8$ of a given degree, that component must be defined 
over $\F_8$.

\subsection{Analysis of cases}

We saw in the previous section that any ``bad'' curve with $27$ points 
must be geometrically reducible.  We therefore organize our discussion 
around the possible lists of degrees for the geometric components of our 
bad curve.

At the outset, we remark that the quadric surface $X^2 + XY + Y^2 + ZW = 
0$ contains no $\F_8$-lines.  Moreover, every $\F_8$-line on the cone $XY 
+ Z^2 = 0$ passes through the vertex of the cone, and in our 
computations we have specifically excluded the cubic surfaces which 
contain the vertex of the cone.  Therefore, every case in which the 
bad curve contains an $\F_8$-line can arise only when the quadric surface 
under consideration is $XY + ZW = 0$, which is isomorphic to 
$\PP^1 \times 
\PP^1$.    We recall that a curve of bidegree $(a,b)$ in $\PP^1 \times 
\PP^1$ 
has arithmetic genus $(a-1)(b-1)$ and intersects a curve of bidegree 
$(c,d)$ exactly $ad + bc$ times.

\noindent \textbf{Degrees $(5,1)$.} By the note at the end of
the preceeding subsection, both components are defined over $\F_8$, and 
so this case can only be found on $XY+ZW=0$. 
By the argument in \ref{red}, if the component of degree $5$ is 
singular, it has at most $15$ points; and since the component of degree 
$1$ has only $9$ points, this is too few points.  Thus the component of 
degree $5$ is nonsingular.   By
lemma \ref{max}, the component of degree $5$ has at most $18$ points, so 
to get a total of $27$ $\F_8$-points must have exactly $18$ points. 

A genus 
$2$ curve of degree $5$ over $\F_8$ with $18$ points has ``defect 
1'' in the terminology of \cite{lauter2}, and the negatives of the 
Frobenius traces are either $5,4$ or $9/2 \pm \sqrt{5}/2$.  By 
criterion
(2.3) of \cite{lauter2}, 
the former cannot occur.  In the second case, one 
checks from the Weil conjectures that the number of $\F_{64}$-points of 
the 
curve is exactly $54$.    

An $\F_8$-line has $9$ points, so the linear 
component and the component
of degree $5$ do not meet over $\F_8$.  Since the components have
bidegrees $(3,2)$ and $(0,1)$, and therefore intersect exactly $3$ times
over the algebraic closure, the two components cannot intersect over
$\F_{64}$ either.  Consequently, in this case we should find exactly $65 + 
54 
= 119$ $\F_{64}$-points on the bad curve.

\noindent \textbf{Degrees $(4,2)$.}  By lemma \ref{max}, there could be at 
most $14+9 = 23$ points on these components, so this case does not occur.

\noindent \textbf{Degrees $(3,3)$.}  If the two components are defined 
over 
$\F_8$, they have at most $9$ $\F_8$-points by lemma \ref{max}; if they 
are not defined over $\F_8$, we draw the same conclusion from lemma 
\ref{field}.  Either way, there are at most $18$ points on these 
components, and so this case does not occur.

\noindent \textbf{Degrees $(4,1,1)$.}  The two lines must be defined over 
$\F_8$, or else there are at most $14+1+1$ points, so we may restrict 
attention to the quadric $XY+ZW=0$.  We note that the 
list of bidegrees must either be $(3,1),(0,1),(0,1)$, or 
$(2,2),(0,1),(1,0)$.  

In the former case, all three components have arithmetic genus $0$, so
have $9$ $\F_8$-points and must not intersect over $\F_8$.  The components 
of bidegree $(3,1)$ and $(0,1)$ intersect three times over the algebraic 
closure, so if they do not intersect over $\F_8$ then they cannot 
intersect over $\F_{64}$.  Since two lines of bidegree $(0,1)$ never 
intersect, there must be a total of $3 \cdot 65 = 195$ points of 
intersection over $\F_{64}$.

In the latter case, the component of bidegree $(2,2)$ has arithmetic 
genus $1$.  The lines of bidegree $(0,1)$ and $(1,0)$ intersect once, and 
so 
have exactly $17$ $\F_8$-points between them.  Thus the curve of genus $1$
must have at least $10$ $\F_8$-points.  We consider each possibility in 
turn, recalling that a singular curve of arithmetic genus $1$ has at most 
$10$ $\F_8$-points.  Note that by Honda-Tate theory, a curve of genus $1$ 
over $\F_8$ does not have $11$ points.  (See Theorem 
4.1 of \cite{waterhouse}.)

\begin{itemize}
\item {\textbf{If the curve of genus $1$ has $10$ points and is 
non-singular,}} then it has $80$ points over $\F_{64}$.  It does not meet 
either line over $\F_8$, but must meet them each in a pair of conjugate 
points over $\F_{64}$.  These intersection points are different for each 
line, as the two lines are distinct.  Since the lines intersect once, the 
total number of $\F_{64}$-points must be $65 + 65 + 80 - 5 = 205$.

\item{\textbf{If the curve of genus $1$ has $10$ points and is singular,}} 
then it has $64$ points over $\F_{64}$.  The rest of our analysis in the 
previous case remains the same, and so the total number of 
$\F_{64}$-points 
must be $65 + 65 + 64 - 5 = 189$.  

\item{\textbf{If the curve of genus $1$ has $12$ points,}} then it has 
$72$ 
points over $\F_{64}$.  The elliptic curve must have two points of 
intersection with the lines over $\F_8$, and so depending on the 
intersection geometry may have either $2$ or $4$ points of intersection 
with the lines over $\F_{64}$.  The total number of $\F_{64}$ points is 
either 
$129 + 72 - 2 = 199$ (if the elliptic curve intersects each line at a 
double-point, or else has a double-point with one line at the intersection 
of the two lines and meets the other line singly there and at one other 
point)  or $129 + 72 
- 4 = 197$ (if the elliptic curve intersects both lines in two distinct 
$\F_8$-points).

\item{\textbf{If the curve of genus $1$ has $13$ points,}} then it has 
$65$ 
points over $\F_{64}$ and must intersect the two lines in three points 
over 
$\F_8$.  The only way this is possible is to pass through the point of 
intersection of the two lines, and to meet each line once more over 
$\F_8$.  Then the total number of points over $\F_{64}$ is $129 + 65 - 3 = 
191$.

\item{\textbf{If the curve of genus $1$ has $14$ points,}} then it has 
$56$ 
points over $\F_{64}$ and has four distinct points of $\F_8$-intersection 
with the lines.  The total number of points over $\F_{64}$ is then $129 + 
56 
- 4 = 181$.
\end{itemize}

\noindent \textbf{Degrees $(3,2,1)$.}  All must be defined over $\F_8$, 
and 
so can occur only in the $XY + ZW=0$ case.  
Each component would have $9$ points, but the component of bidegree 
$(1,1)$ must meet the linear component, so we cannot reach as many as $27$ 
$\F_8$-points.

\noindent \textbf{Degrees $(2,2,2)$.}  All must be defined over $\F_8$, or 
else we have at most $9 + 4 + 4 < 27$ $\F_8$-points.  Each component has 
$9$ points, and is the intersection of a plane with our quadric.  Hence 
any two of the components intersect in $2$ points over $\F_{64}$, and so 
have $128$ $\F_{64}$-points between them.  The third component has either 
$65-2$ or $65-4$ points not on either of the first two, and so there are 
either $189$ or $191$ $\F_{64}$-points in total.  Note that this is the 
only 
case in which we are not limited to the split non-degenerate quadric.

\noindent \textbf{Degrees $(3,1,1,1)$.} If the lines are not all defined 
over 
$\F_8$, then there are at most $9+9+1+1<27$ points.  The bidegrees must be 
$(1,2),(1,0),(1,0),(0,1)$ and so there are at most $7$ points of 
intersection between the components.  Then there are at least $36-7 > 27$  
$\F_8$-points, which is too many, and so this case cannot occur.

\noindent \textbf{Degrees $(2,2,1,1)$.} Again, every component must be 
defined over $\F_8$, and the bidegrees are $(1,1),(1,1),(0,1),(1,0)$.  
Once again there are too many points.

\noindent \textbf{Degrees $(2,1,1,1,1)$.} At least two of the lines must 
be 
defined over $\F_8$.  So, if not all of the lines are defined over $\F_8$, 
then precisely two are not.  The two lines not defined over $\F_8$ 
would either both have bidegree $(1,0)$ or both have bidegree $(0,1)$, and 
so  
would not 
meet; therefore they could not contain any $\F_8$-points, as the lines are 
Galois-conjugate and any $\F_8$-points on them would lie in their 
intersection.  Therefore, since the curve of bidegree $(1,1)$ intersects 
the two $\F_8$-lines, the configuration could contain at most $27-2=25$ 
$\F_8$-points.  On the other hand, if all the lines are defined over 
$\F_8$, there are far too many $\F_7$-points.  So this case cannot occur.

\noindent \textbf{Degrees $(1,1,1,1,1,1)$.}  If four of the lines are 
defined over $\F_8$, then there are too many points; and if there are 
only two, then there are too few points.  However, it is possible that 
the three lines of (say) bidegree $(1,0)$ could be defined over $\F_8$, 
while the three lines of bidegree $(0,1)$ could  be defined over 
$\F_{512}$.  
Then over $F_{64}$ there would be exactly $3 \cdot 65 = 195$ points.

To summarize: on any of the quadrics, our bad curve may decompose into 
three plane quadric curves over $\F_8$.  In this case there are either 
$189$ or $191$ $\F_{64}$-points on the bad curve.  This is the only 
possibility on the degenerate and non-split non-degenerate quadrics.  In 
the split non-degenerate case, we have the following additional 
possibilities:

\begin{itemize}
\item The bad curve has two components, both defined over $\F_8$, one of 
bidegree $(3,2)$ and 
one of bidegree $(0,1)$.  In this case there are $119$ $\F_{64}$-points.

\item The bad curve has three components, all defined over 
$\F_8$, one of bidegree $(3,1)$ 
and two lines of bidegree $(0,1)$.  In this case there are $195$ 
$\F_{64}$-points.

\item The bad curve has three components, all defined over $\F_8$, one of 
bidegree $(2,2)$ 
and lines of bidegree $(0,1)$ and $(1,0)$.  In this case, there are $189$, 
$205$, $199$, $197$, $191$, or $181$ $\F_{64}$-points, depending on either 
the curve of bidegree $(2,2)$ is singular with $10$ $\F_8$-points, or 
non-singular with $10$, $12$, $12$, $13$, or $14$ $\F_8$-points 
respectively.

\item The bad curve has six linear components, three defined over $\F_8$ 
and three defined over $\F_{512}$.  In this case there are $195$ 
$\F_{64}$-points.
\end{itemize}

\subsection{Examples}

Scouring our computer calculations, 
we have found an example of each of the possibilities for bad curves
enumerated in the
previous section, and so all of these possibilities do indeed occur.  We
give a few of these examples here; the interested reader may refer to
math.NT/0201226 at {\sf http://arXiv.org} or to
\begin{center}
{\sf http://www.math.mcgill.ca/$\sim$dsavitt/curves/examples.dvi} 
\end{center}
for the full list.  (This file is also available in {\sf .ps} and {\sf .pdf}
format.)

Recall that $\eta \in \F_8$ is a 
chosen root of $\eta^3 + \eta + 1 = 0$.  Let $\beta$ be a generator of 
$\F_{64}^{\times}$ such that $\beta^9 = \eta$.  Each intersection described 
below has exactly $27$ points over $F_8$.

\begin{itemize} 

\item The intersection of $XY+ZW=0$ with the cubic 
$X^2W + \eta XYW + \eta^{-1} XZW + \eta^{-3} XW^2 + \eta Y^2Z + Y^2W + 
\eta^{-2} YZ^2 + \eta^{-1} YZW + YW^2 = 0$
contains the line 
$[X:0:Z:0]$ and a component of degree $5$, and has $119$ points over 
$\F_{64}$.

\vskip 0.1cm

\item The intersection of $XY+ZW=0$ with the cubic
$$(\eta Y + Z)(YZ + XZ + \eta XW + \eta^{-1} W^2 + \eta ZW + \eta^{-1} YW)$$ 
contains the lines 
$[\eta 
W : Y : \eta Y : W]$ and $[X : 0 : 0 : W]$. 
The intersection of $XY+ZW=0$ with $YZ + XZ + \eta XW + 
\eta^{-1} W^2 + \eta ZW + \eta^{-1} YW = 0$ is an elliptic curve with $12$ 
$\F_8$-points and $72$ $\F_{64}$-points.  It meets 
the line $[X: 0 : 0 : W]$ at the two 
points $[1:0:0:0]$ and $[1:0:0:\eta^2]$, and meets the line $[\eta W : Y : 
\eta Y : W]$ at the two Galois-conjugate points $[ \beta^{59} : 1 : 
\beta^9 : \beta^{50}]$ 
and $[ \beta^{31} : 1 : \beta^9 : \beta^{22} ]$.  The 
intersection of the quadric 
and the cubic has $197$ points 
over $\F_{64}$.

\vskip 0.1cm

\item The intersection of $XY+ZW=0$ with the cubic $\eta^{-2} X^2Z + 
\eta^3 XYZ + \eta^3 XYW + \eta^{-2} XZ^2 + \eta^3 XZW + Y^2W + \eta^3 YZW + 
YW^2 = 0$ contains the three non-intersecting lines $[0:Y:Z:0]$, 
$[X:0:0:W]$, and $[X:Y:X:Y]$ and three lines defined over $\F_{512}$.  The 
intersection has $195$ points over $\F_{64}$.

\end{itemize}

\section*{{\bf APPENDIX}}

\begin{center}
by Kristin Lauter, Microsoft Research
\end{center}

\markboth{\textsc{Kristin Lauter}}{\textsc{The maximum number of points on 
a 
curve of genus $4$ over $\F_8$ is $25$}}
\pagestyle{myheadings}

\subsection*{A.1. Introduction}

The purpose of this appendix is to give a list of the possible zeta
functions for curves with defect $3$. As a special case, we will show
that there is no genus $4$ curve over $\F_8$ with
$26$ rational points.

\subsection*{A.2. Definitions}

By a {\it curve} over $\F_q$, we mean a smooth, projective,
absolutely irreducible curve.  For a curve, $C$,
let $g=g(C)$ denote the genus, and $N(C)$ denote the
number of rational points over $\F_q$.
A curve $C$ has {\it defect $k$} if
it fails to meet the Serre-Weil bound by exactly $k$:
$$N(C)=q+1+gm-k,$$
where $$m=[2\sqrt{q}].$$
The {\it zeta function} of a curve over $\F_q$ is defined as a power 
series,
but it is known that it is a rational function, and can be
written in the form
$$\frac{h(t)}{(1-t)(1-qt)},$$ where
$$h(t) = \prod_{i=1}^g (1-\alpha_it)(1- \bar{\alpha_i}t)$$
is a polynomial with coefficients in $\Z$, and $\alpha_i$ and 
$\bar{\alpha_i}$
are algebraic integers with complex absolute value $\sqrt{q}$.
We say that a curve has {\it zeta function of type $(x_1,...x_g)$} if
$x_i = -(\alpha_i + \bar{\alpha_i}), \quad i=1,...,g.$
Define the polynomial $P(t)$:
$$P(t) = \prod_{i=1}^{g}(t-(m+1-x_i)),$$
and the set $F_k$:
$$F_k = \{ t^d + a_1 t^{d-1} + ... + a_d \in \Z[t] \mid - a_1 = d+k,
\; {\rm all \:  roots \: positive \: reals} \}.$$
The $m+1-x_i$ are totally positive algebraic integers,
so if $$\sum_{i=1}^{g} x_i = gm-k,$$
then $P(t) \in F_k,$ since $\deg P =g$, and  $- a_1 = g+k$.
We say $P(t)$ is a {\it polynomial of defect $k$}.

\subsection*{A.3. Defect 3}

Using the method of Smyth as explained in \cite{Serre-Harvard} or
Section 2 of \cite{LS}, we restrict the possibilities for the type of the
zeta function
for defect $3$ curves by making a list of the possibilities
for the irreducible factors of the polynomials $P(t)$.

The possibilities are divided into four types given in the
following four tables.

\begin{itemize}

\item {\bf Type 1} is an irreducible polynomial of defect $3$ and the rest 
of the
factors are made up of defect $0$ polynomials.
For $k=0$, the defect $0$ polynomial is $P(t)=(t-1)$, so the $x_i$
corresponding to this factor is $x_i = m$.

\item {\bf Type 2} is an irreducible polynomial of defect $2$ combined 
with the
defect $1$ polynomial $(t-2)$
and copies of the defect $0$ polynomial $(t-1)$.

\item {\bf Type 3} is an irreducible polynomial of defect $2$ combined
with the defect $1$ polynomial $(t^2-3t+1)$
and copies of the defect $0$ polynomial $(t-1)$.

\item {\bf Type 4} consists of the four possible combinations of the two 
defect $1$
polynomials with the rest of the factors equal to the defect $0$
polynomial $(t-1)$.
\end{itemize}

\begin{table}[ht]
\begin{center}

\caption{Possibilities for $P(t)$ and $(x_1,...,x_g)$ for defect $3$: Type 
1}

\begin{tabular}{|c|c|c|c|c|c|c|}
\hline
\# & deg & coefficients & $(x_1,...,x_g)$ & $g\ge\:?$ & 
$\{2\sqrt{q}\}\ge\:?$\\
\hline
1. & 4 & 1 -7 14 -8 1 & & $g \ge 4$& 0.827\dots  \\
\hline
2. & 4 & 1 -7 13 -7 1 & & $g \ge 4$& 0.772\dots  \\
\hline
3. & 3 & 1 -6 5 -1 & & $g \ge 3$ & 0.692\dots  \\
\hline
4. & 3 & 1 -6 7 -1 & & $g \ge 3$ & 0.834\dots  \\
\hline
5. & 3 & 1 -6 8 -1 & & $g \ge 3$ & 0.860\dots  \\
\hline
6. & 3 & 1 -6 8 -2 & & $g \ge 3$ & 0.675\dots  \\
\hline
7. & 3 & 1 -6 9 -1 & & $g \ge 3$ & 0.879\dots  \\
\hline
8. & 3 & 1 -6 9 -3 & & $g \ge 3$ & 0.532\dots  \\
\hline
9. & 2 & 1 -5 5 & ($m,\dots,m-\frac{3\pm\sqrt{5}}{2}$) & $g \ge 2$ &   \\
\hline
10.& 2 & 1 -5 3&($m,\dots,m-\frac{3\pm\sqrt{13}}{2}$) &$g\ge 2$ & 
0.302\dots \\
\hline
11. & 2 & 1 -5 2 & ($m,\dots,m-\frac{3\pm\sqrt{17}}{2}$)&$g\ge 
2$&0.561\dots \\
\hline
12. & 2 & 1 -5 1 & ($m,\dots,m-\frac{3\pm\sqrt{21}}{2}$)&$g \ge 
2$&0.791\dots\\
\hline
13. & 1 & 1 -4 & ($m,\dots,m-3$) & $g \ge 1$ & 0 \\
\hline

\end{tabular}
\end{center}
\end{table}

\begin{table}[ht]
\begin{center}

\caption{Possibilities for $P(t)$ and $(x_1,...,x_g)$ for defect $3$: Type 
2}

\begin{tabular}{|c|c|c|c|c|c|}
\hline
\# & deg & coefficients & $(x_1,...,x_g)$ & $g\ge\:?$ & 
$\{2\sqrt{q}\}\ge\:?$\\
\hline
14. & 3 & 1 -5 6 -1 & & $g \ge 4$& 0.8019\dots  \\
\hline
15. & 2 & 1 -4 2 &($m-(1\pm\sqrt{2}),m-1,m,\dots$)&$g \ge 3$& 0.414\dots 
\\
\hline
16. & 2 & 1 -4 1 &($m-(1\pm\sqrt{3}),m-1,m,\dots$)  & $g \ge 3$& 
0.732\dots \\
\hline
17. & 1 & 1 -3 & ($m-2,m-1,m,\dots$) & $g \ge 2$ & 0 \\
\hline

\end{tabular}
\end{center}
\end{table}

\begin{table}[ht]
\begin{center}

\caption{Possibilities for $P(t)$ and $(x_1,...,x_g)$ for defect $3$: Type 
3}

\begin{tabular}{|c|c|c|c|c|c|}
\hline
\# & deg & coefficients & $(x_1,...,x_g)$ & $g\ge\:?$ & 
$\{2\sqrt{q}\}\ge\:?$\\
\hline
18. & 3 & 1 -5 6 -1 & & $g \ge 5$& 0.8019\dots  \\
\hline
19. & 2 & 1 -4 2 &($m-(1\pm\sqrt{2}),m-\frac{1\pm\sqrt{5}}{2},m,\dots$)
& $g \ge 4$& 0.618\dots  \\
\hline
20. & 2 & 1 -4 1 &($m-(1\pm\sqrt{3}),m-\frac{1\pm\sqrt{5}}{2},m,\dots$)
& $g \ge 4$& 0.732\dots  \\
\hline
21. & 1& 1 -3 & ($m-2,m-\frac{1\pm\sqrt{5}}{2},m,\dots$) & $g \ge 3$
& 0.618\dots \\
\hline

\end{tabular}
\end{center}
\end{table}

\begin{table}[ht]
\begin{center}

\caption{Possibilities $(x_1,...,x_g)$ for defect $3$: Type 4}

\begin{tabular}{|c|c|c|c|}
\hline
\# & $(x_1,...,x_g)$ & $g\ge\:?$ & $\{2\sqrt{q}\}\ge\:?$\\
\hline
22. &($m-1,m-1,m-1,m,\dots$) & $g \ge 3$& 0  \\
\hline
23. &($m-\frac{1\pm\sqrt{5}}{2},m-1,m-1,m\dots$) & $g \ge 4$& 0.618\dots  
\\
\hline
24. & ($m-\frac{1\pm\sqrt{5}}{2},m-\frac{1\pm\sqrt{5}}{2},m-1,m,\dots$)
& $g \ge 5$ & 0.618\dots \\
\hline
25. &($m-\frac{1\pm\sqrt{5}}{2},m-\frac{1\pm\sqrt{5}}{2},
m-\frac{1\pm\sqrt{5}}{2},m,\dots$)  & $g \ge 6$& 0.618\dots  \\
\hline

\end{tabular}
\end{center}
\end{table}

For each pair $(q,g)$,there could be a number of reasons why
an entry in the above tables does not correspond to the zeta
function of  a curve.

Using the following three reasons from Section 2 of \cite{LS}
we can eliminate many of the entries from the tables.

\begin{list}{}{}
\item (2.1)  The absolute value of each $x_i$ must be less than 
$2\sqrt{q}$.
\item (2.2)  The number of places of degree $d$ on a curve is 
non-negative.
\item (2.3)  The numerator of the zeta function of a curve is not 
decomposable.
\end{list}

The last column in each table indicates the restriction
that comes from reason (2.1):
$\{2\sqrt{q}\} \ge 1-x$, where $x$ is the smallest root of $P(t)$.


\noindent \textbf{Proposition A.1} \textit{
The following entries from the tables do 
not 
correspond
to the zeta function of a curve for reason (2.3).
\begin{itemize}{}{}
\item \#17 for genus $g \ge 2$,
\item \#9,10,21 for genus $g \ge 3$,
\item \#3,4,6,8,14,15,19,20,22,23 for genus $g \ge 4$,
\item \#1,2,18,24 for genus $g \ge 5$,
\item \#25 for genus $g \ge 7$.
\end{itemize}
}

\begin{proof} 
 For each entry, it suffices to factor the corresponding polynomial
$$F(T)=\prod_{i=1}^g (T-(\alpha_i+\bar{\alpha_i}))=\prod_{i=1}^g (T+x_i)$$
into two factors, $f(T)$ and $g(T)$ such that the resultant of $f$ and
$g$ is $\pm 1$ (see Lemma 4.1, \cite{L}).  For example, for entry \#8,
the resultant of
$$T^3+(3m-3)T^2+(3m^2-6m)T+m^3-3m^2+1$$
and $(T+m)$ is $-1$, so entry \#8 is not possible for $g \ge 4$.
For entry \#19,
$${\rm resultant}(T^2+(2m-2)T+m^2-2m-1, T^2+(2m-1)T+m^2-m-1) = -1,$$
so this entry is not possible for $g = 4$, and
$${\rm resultant}((T^2+(2m-2)T+m^2-2m-1)(T+m), T^2+(2m-1)T+m^2-m-1) = 1,$$
so it is not possible for $g > 4$ either.  The decomposition of other 
entries
is similar.
\end{proof}

\noindent \textbf{Proposition A.2} \textit{
Entry \#11 does not correspond
to the zeta function of a curve for
$$g > \frac{q^2-q+8m^2-10m-16}{5m^2-7m-2q}$$
for reason (2.2).
}

\begin{proof} The proof is similar to the proof of Proposition 1 in 
\cite{LS}.
The coefficients of the polynomial
$$(T+m-\frac{3+\sqrt{17}}{2})(T+m-\frac{3-\sqrt{17}}{2})(T+m)^{g-2}$$
can be computed in two ways: as binomial coefficients
or via Newton's relations between the elementary symmetric
functions, $\{b_n\}$,  and the power functions,
$$s_n = \sum_{i=1}^g (\alpha_i + \bar{\alpha_i})^n.$$
Using the identity
$$b_2 = \frac{1}{2}(s_1^2 -s_2),$$
and equating the coefficients of the $g-2$ term computed
in the two ways yields:
\begin{multline*}
\frac{(g-2)(g-3)}{2}m^2+(g-2)m+(m^2-3m-2) \\
= \frac{1}{2}((gm-3)^2-(q^2+1-(q+1+gm-3+2a_2)+2gq)) \,,
\end{multline*}
where $a_2$ is the number of places of degree $2$ on the curve.
By reason (2.2), we must have $a_2 \ge 0$, so rearranging yields
the desired inequality.
\end{proof}

\noindent \textbf{Proposition A.3} \textit{
Entry \#13 does not correspond
to the zeta function of a curve for
$$g > \frac{q^2-q+6m-6}{m^2+m-2q}$$
for reason (2.2).  In general, $(m,m,\dots,m-k)$ does not correspond
to the zeta function of a defect $k$ curve for
$$g > \frac{q^2-q+2km+k-k^2}{m^2+m-2q}.$$
}

\begin{proof} The proof is similar to the proof of proposition A.2 above.
\end{proof}

\noindent \textbf{Remark.}  Similar bounds on the genus can be obtained 
for entries \#5,7,12,16.

\noindent \textbf{Proposition A.4} \textit{
If $q$ is an even power of a prime, then the 
only defect $3$ 
curves with genus
$g>3$ have zeta function of type $(m,...,m,m-3)$.  For $g=3$,
$(m-1,m-1,m-1)$ is possible in some cases. For
$$g> \frac{q^2-q+6m-6}{m^2+m-2q},$$ defect $3$ is not possible.
}

\begin{proof} This follows from reason (2.1) and the fact that entries 
\#17
and \#22 are impossible by reason (2.3) for $g \ge 2$ and $g \ge 4$
respectively.  The last statement then follows from proposition A.3.
\end{proof}

\noindent \textbf{Theorem A.5} \textit{
\label{t26} There does not exist a genus $4$ curve over $\F_8$ 
with $26$ 
$\F_8$-points.
}

\begin{proof}  When $q=8$, $$\{2\sqrt{q}\} \approx 0.6568,$$ so using the
above tables, we see that the only zeta function types possible after 
applying
proposition A.1 are:  \#11 and \#13.  By proposition A.2, \#11 is not 
possible
since $g=4 > \frac{95}{37}$.  For \#13, the bound on $g$ from proposition 
A.3
is $\frac{40}{7}>4$, but \#13 is not possible for a different reason in 
this
case.  Here $q=2^3$ and $m=5$, so $m-3=2$. By Honda-Tate theory, when
$q=p^e$ is an odd power of a prime, the only possible
values for the trace of an elliptic curve which are
divisible by the characteristic are: (see \cite{Water}, p.536)
\begin{list}{}{}
\item $0$, for all $p$, or
\item $p^{\frac{e+1}{2}}$, for $p=2$ or $p=3$.
\end{list}{}{}
Since an elliptic curve with trace $2$ does not exist over $\F_8$,
an abelian variety over $\F_8$ of type $(5,5,5,2)$ does not exist either.
\end{proof}

Theorem A.5 was presented at the Journ\'ees Arithm\'etiques in Rome 
in July, 1999, and at the Arizona Winter School in March, 2000.

\footnotesize

\small

\noindent \textsc{Department of Mathematics, McGill University, and CICMA}

\vskip -.2cm

\noindent \textsf{dsavitt@math.mcgill.ca}

\vskip -.1cm

\noindent \textsc{Microsoft Research}

\vskip -.2cm
\noindent \textsf{klauter@microsoft.com}


\begin{thebibliography}{[99]}

\bibitem[ACGH]{ACGH} Arbarello, E., M. Cornalba, P.A. Griffiths, and
J. Harris. {\em Geometry of Algebraic Curves, Volume I.}  New York:
Springer-Verlag, 1985.

\bibitem[Arf]{arf} Arf, Cahit.  {\em Untersuchungen \"uber quadratische 
Formen in K\"orpern der Charakteristik 2.  I.}  J. Reine Angew. Math. 
\textbf{183} (1941), pp. 148--167.

\bibitem[BCP]{bcp} Bosma, W., J. Cannon, and C. Playoust. {\em The Magma 
algebra system. I. The user language.}  Computational algebra and number 
theory (London, 1993).  J. Symbolic Comput. \textbf{24} (1997) 3-4, pp. 
235--265.

\bibitem[GV]{gv} van der Geer, Gerard and Marcel van der Vlugt. {\em 
Tables 
of curves with many points.}  Available at: 
{\sf http://www.science.uva.nl/\~{}geer/tables-mathcomp9.ps}

\bibitem[Har]{hartshorne}  Hartshorne, Robin.  {\em Algebraic Geometry.}  
New York: Springer-Verlag, 1977.

\bibitem[HS]{sh} Hindry, Marc and Joseph H. Silverman.  {\em Diophantine 
Geometry, An Introduction.}  New York: Springer-Verlag, 2000.

\bibitem[Lau1]{lauter} Lauter, Kristin.  {\em Improved upper bounds for 
the 
number of rational points on algebraic curves over finite fields.} C.R. 
Acad. Sci. Paris \textbf{328}, S\'erie I (1999), pp. 1181--1185.

\bibitem[Lau2]{lauter2} Lauter, Kristin.  {\em Geometric methods for 
improving the upper bounds on the number of rational points on algebraic 
curves over finite fields.} With an appendix in French by J.-P. Serre.  J. 
Algebraic Geom. \textbf{10} (2001), no. 1, pp. 19--36.

\bibitem[Se1]{serre} Serre, J-P. {\em Algebraic Groups and Class Fields.}  
New York: Springer-Verlag, 1988.

\bibitem[Se2]{serre2} Serre, J-P. {\em Nombre de points des courbes
alg\'ebriques sur ${\mathbb F}_q$.}  S\'em. de Th\'eorie des Nombres de
Bordeaux, 1982/83, exp. no. 22.  (= Oeuvres III, No. 129, p. 664-668).

\bibitem[Wat]{waterhouse} Waterhouse, William C. {\em Abelian Varieties 
over Finite Fields.}  Ann. scient. \'Ec. Norm. Sup., $4^{e}$ s\'erie, t.2, 
1989, pp. 521--560.

\end{thebibliography}

\begin{thebibliography}{99}

\bibitem{L} K. Lauter, {Non-existence of a curve over ${F}\sb 3$ of genus
$5$ with $14$ rational points}, {\it Proc. Amer. Math. Soc.} {\bf 128},
369-374 (2000).

\bibitem{LS} K. Lauter, with an Appendix by J-P. Serre,
{\em Geometric Methods for Improving the Upper Bounds on the Number of
Rational Points on Algebraic Curves over Finite Fields},
Journal of Algebraic Geometry {\bf 10} (2001), no. 1, 19-36

\bibitem{Serre-Harvard} J.-P. Serre, {\em Rational Points on Curves over
Finite Fields}. Notes by F. Gouvea of lectures at Harvard University,
1985.

\bibitem{Smyth} C. Smyth, {\em Totally Positive Algebraic Integers
of Small Trace}, Ann. Inst. Fourier, Grenoble {\bf 33}, 3 (1984), 1-28.

\bibitem{Water} W.C. Waterhouse, {\em Abelian Varieties over Finite
Fields}, Ann. scient. \'Ec. Norm. Sup., $4^e$ s\'erie, t. 2, 1969, p. 
521-560.

\end{thebibliography}
\end{document}